\documentclass{amsart}
\usepackage[top=1 in,bottom=1.25 in,left=2cm, right=2cm]{geometry}
\usepackage{amsmath, amsthm, amscd, amsfonts, amssymb, color, hyperref}
\usepackage[dvips]{graphicx}
\usepackage{tikz}
\usepackage{float}
\usepackage{perpage}
\usepackage[T1]{fontenc}
\newtheorem{thm}{Theorem}[section]

\theoremstyle{definition}

\numberwithin{equation}{section}
\begin{document}

\title[cycles, paths and certain trees]{A Conjecture about spectral distances between cycles, paths and certain trees}%
\author[A. Abdollahi]{Alireza Abdollahi}%
\address{Department of Pure Mathematics, Faculty of Mathematics and Statistics, University of Isfahan, Isfahan 81746-73441, Iran; and School of Mathematics, Institute for Research in Fundamental Sciences (IPM), P.O. Box 19395-5746, Tehran, Iran}%
\email{a.abdollahi@math.ui.ac.ir}%
\author[N. Zakeri]{Niloufar Zakeri}
\address{Department of Pure Mathematics, Faculty of Mathematics and Statistics, University of Isfahan, Isfahan 81746-73441, Iran}%
\email{zakeri@sci.ui.ac.ir}%
\thanks{}%
\subjclass[2010]{05C50; 05C31}%
\keywords{Adjacency spectra of graphs; Spectral distance of graphs; Path; Cycle}%

%\date{}%
%\dedicatory{}%
%\commby{}%
% ----------------------------------------------------------------
\begin{abstract}
We confirm the following conjecture which has been proposed in  [{\em Linear Algebra and its Applications}, {\bf 436} (2012), No. 5,
1425-1435.]:
$$  0.945\approx\displaystyle\lim_{n\longrightarrow \infty}\sigma(P_n,Z_n)=\displaystyle\lim_{n\longrightarrow \infty}\sigma(W_n,Z_n)=\frac{1}{2}\displaystyle\lim_{n\longrightarrow \infty}\sigma(P_n,W_n);\ \displaystyle\lim_{n\longrightarrow \infty}\sigma(C_{2n},Z_{2n})=2,$$
where $\sigma(G_1,G_2)=\sum_{i=1}^n |\lambda_i(G_1)-\lambda_i(G_2)|$ is the spectral distance between $n$ vertex non-isomorphic graphs $G_1$ and $G_2$ with adjacency spectra $\lambda_1(G_i) \geq \lambda_2(G_i) \geq \cdots \geq \lambda_n(G_i)$ for $i=1,2$, and  $P_n$ and $C_n$ denote the path and cycle on $n$ vertices, respectively; $Z_n$ denotes the coalescence of $P_{n-2}$ and $P_3$ on  one of the  vertices of degree 1 of $P_{n-2}$ and the vertex of degree $2$ of $P_3$; and $W_n$  denotes the coalescence of $Z_{n-2}$ and $P_3$ on   the  vertex of degree 1 of $Z_{n-2}$ which is adjacent to a vertex of degree $2$ and the vertex of degree $2$ of $P_3$.
\end{abstract}
\maketitle
% ----------------------------------------------------------------
\section{\bf Introduction and Results}
All graphs considered here are simple, that is finite and undirected without loops and multiple edges. By the eigenvalues of $G$, we mean those of its adjacency matrix.
% We denote by $Spec(G)$ the multiset of the eigenvalues of the graph $G$.\\

%To find some applications of the cospectrality of graphs, we refer to \cite{clss,wz,zp}.

In \cite{js}, the following conjecture has been proposed:
$$ 0.945\approx\displaystyle\lim_{n\longrightarrow \infty}\sigma(P_n,Z_n)=\displaystyle\lim_{n\longrightarrow \infty}\sigma(W_n,Z_n)=\frac{1}{2}\displaystyle\lim_{n\longrightarrow \infty}\sigma(P_n,W_n);\ \displaystyle\lim_{n\longrightarrow \infty}\sigma(C_{2n},Z_{2n})=2,$$
%It has been considered as an extension of the results given in Theorem 3.5 in \cite{js}.
where $\sigma(G_1,G_2)=\sum_{i=1}^n |\lambda_i(G_1)-\lambda_i(G_2)|$ is the spectral distance between $n$ vertex non-isomorphic graphs $G_1$ and $G_2$ with adjacency spectra $\lambda_1(G_i) \geq \lambda_2(G_i) \geq \cdots \geq \lambda_n(G_i)$ for $i=1,2$. 
In this paper we completely prove the conjecture.\\
 The spectral distance between graphs is studied based on different matrix representations of graphs and norms in many papers  for example see \cite{ajo,ds,jzs,lld}.

Let us first introduce some notations.

The coalescence of two graphs $G$ and $H$ is obtained  by identifying a vertex $u$ of $G$ with a vertex $v$ of $H$; $Z_n$ denotes the coalescence of $P_{n-2}$ and $P_3$ on  one of the  vertices of degree 1 of $P_{n-2}$ and the vertex of degree $2$ of $P_3$; and $W_n$  denotes the coalescence of $Z_{n-2}$ and $P_3$ on   the  vertex of degree 1 of $Z_{n-2}$ which is adjacent to a vertex of degree $2$ and the vertex of degree $2$ of $P_3$. They are depicted in Figure 1 (Their spectra can be found in \cite{cds}, p. 77).

Our main result is as follows.

\begin{thm}\label{cz}
\
\begin{enumerate}
\item $\displaystyle\lim_{n\longrightarrow \infty}\sigma(C_{2n},Z_{2n})=2.$\\

\item $\displaystyle\lim_{n\longrightarrow \infty}\sigma(P_n,Z_n)=\frac{8-8\sqrt{2}+2\pi}{\pi}\approx 0.945.$\label{pz}\\

\item $\displaystyle\lim_{n\longrightarrow \infty}\sigma(W_n,Z_n)=\frac{8-8\sqrt{2}+2\pi}{\pi}\approx 0.945.$\label{wz}\\

\item $ \displaystyle\lim_{n\longrightarrow \infty}\sigma(P_n,Z_n)=\displaystyle\lim_{n\longrightarrow \infty}\sigma(W_n,Z_n)=\frac{1}{2}\displaystyle\lim_{n\longrightarrow \infty}\sigma(P_n,W_n)=\frac{8-8\sqrt{2}+2\pi}{\pi}\approx 0.945.$  
\end{enumerate}

\end{thm}
\begin{figure}
\begin{center}
\includegraphics{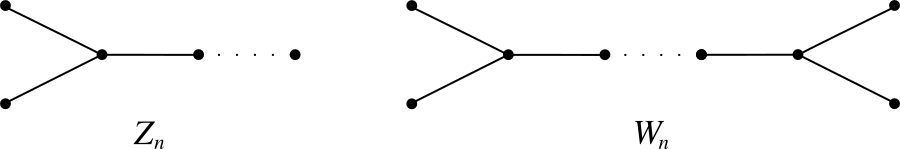}
\caption{Graphs $Z_n$ and $W_n$.}
\label{fig:g}
\end{center}
\end{figure}

\section{\bf Proof of the conjecture}

To prove the following results, we need the Table 1.
 
\begin{center}
\begin{table}

\begin{tabular}{|c|l|}\hline
Graph & Eigenvalues \\ \hline
$P_n$& $2\cos(\frac{k\pi}{n+1}),k=1,\ldots,n$\\
$C_n$& $2\cos(\frac{2k\pi}{n}),k=1,\ldots,n$\\
$Z_n$& $0, 2\cos(\frac{(2k-1)\pi}{2n-2}),k=1,\ldots,n-1$\\
$W_n$& $2,0,0,-2, 2\cos(\frac{k\pi}{n-3}),k=1,\ldots,n-4$\\\hline
\end{tabular}

\caption{Some specific graphs and their spectra}
\end{table}

\end{center}

\noindent{\bf Proof of Theorem \ref{cz}.} 
\begin{enumerate}

\item Let $n\geq 2$. By Table 1,
$$\lambda_1(C_{2n})>\lambda_k(C_{2n})=\lambda_{k+1}(C_{2n})>\lambda_{k+2}(C_{2n})=\lambda_{k+3 }(C_{2n})>\lambda_{2n}(C_{2n}), \ \ k=2,4,6,\ldots,2n-4.$$
Assume that $n$ is odd. We have

$$\lambda_k(C_{2n})>\lambda_k(Z_{2n})>\lambda_{k+1}(Z_{2n})>\lambda_{k+1 }(C_{2n}), \ \ k=1,3,5,\ldots,2n-1.$$
If $n$ is even, then $\lambda_k(C_{2n})=\lambda_k(Z_{2n})$ for $k=n,n+1$ and the other eigenvalues are the same as the case that $n$ is odd.

Since $C_{2n}$ and $Z_{2n}$ are bipartite, by using the above (in)equalities, we have
\begin{eqnarray*}
\sigma(C_{2n},Z_{2n})&=&2(\lambda_1(C_{2n})-\lambda_1(Z_{2n}))+\displaystyle\sum_{k=2}^{2n-1}|\lambda_k(C_{2n})-\lambda_{k}(Z_{2n})|\\
&=&2(\lambda_1(C_{2n})-\lambda_1(Z_{2n}))+\displaystyle\sum_{{k=2}\; {\text{, $k$ is even}\;}}^{2n-2}\lambda_k(Z_{2n})-\lambda_k(C_{2n})+\lambda_{k+1}(C_{2n})-\lambda_{k+1}(Z_{2n})\\
&=&2(\lambda_1(C_{2n})-\lambda_1(Z_{2n}))+2\displaystyle\sum_{k=2}^{n-1}(-1)^k\lambda_k(Z_{2n})\\
&=&2\lambda_1(C_{2n})+2\displaystyle\sum_{k=1}^{n-1}(-1)^k\lambda_k(Z_{2n})\\
&=&4+4\displaystyle\sum_{k=1}^{n-1}(-1)^k\cos(\frac{(2k-1)\pi}{4n-2}).
\end{eqnarray*}
Since 
$\displaystyle\lim_{n\rightarrow \infty}\displaystyle\sum_{k=1}^{n-1}(-1)^k\cos(\frac{(2k-1)\pi}{4n-2})=\frac{-1}{2},$ we have
$$\lim_{n\rightarrow \infty}\sigma(C_{2n},Z_{2n})=2.$$

\item Let $n\geq 4$ and $a_n:=\sigma(P_n,Z_n)$. We prove
$$\displaystyle\lim_{n\longrightarrow \infty}a_{4n+1}=\displaystyle\lim_{n\longrightarrow \infty}a_{4n+2}=\displaystyle\lim_{n\longrightarrow \infty}a_{4n+3}=\displaystyle\lim_{n\longrightarrow \infty}a_{4n}=\frac{8-8\sqrt{2}+2\pi}{\pi},$$
 which follows that $a_n$ is convergent to $\frac{8-8\sqrt{2}+2\pi}{\pi}$. \\
 
 Suppose that $\lambda_1(P_n)\geq\cdots\geq\lambda_n(P_n)$ and $\lambda_1(Z_n)\geq\cdots\geq\lambda_n(Z_n)$ are the eigenvalues of $P_n$ and $Z_n$, respectively. If  $n \equiv 1 \pmod 4 $, then 
$$\lambda_k(Z_n)>\lambda_k(P_n),\ \ \  k=1,\ldots, \frac{n-1}{4},$$
 
 $$\lambda_k(Z_n)<\lambda_k(P_n),\ \ \ k=\frac{n+3}{4},\ldots, \frac{n-1}{2}, $$
 and $$\lambda_\frac{n+1}{2}(Z_n)=\lambda_\frac{n+1}{2}(P_n).$$
 Thus, by Table 1 and the bipartiteness of $P_n$ and $Z_n$, we have
 \begin{eqnarray*}
\sigma(P_n,Z_n)&=&2\big(\displaystyle\sum_{k=1}^{\frac{n-1}{4}}(\lambda_k(Z_n)-\lambda_{k}(P_n))+\displaystyle\sum_{k=\frac{n+3}{4}}^{\frac{n-1}{2}}(\lambda_k(P_n)-\lambda_{k}(Z_n))\big )\\
&=&4\big(\displaystyle\sum_{k=1}^{\frac{n-1}{4}}(\cos(\frac{(2k-1)\pi}{2n-2})-\cos(\frac{k\pi}{n+1}))
+\displaystyle\sum_{k=\frac{n+3}{4}}^{\frac{n-1}{2}}(\cos(\frac{k\pi}{n+1})-\cos(\frac{(2k-1)\pi}{2n-2}))\big).
\end{eqnarray*}
By using a simple code (e.g., see \ref{fig:m}) in symbolic computational software Maple \cite{m}, we have
\begin{figure}
\begin{center}
\includegraphics{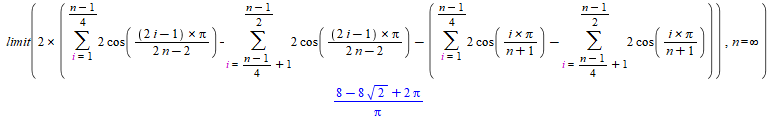}
\caption{The Maple code}
\label{fig:m}
\end{center}
\end{figure}

  $$\displaystyle\lim_{n\longrightarrow \infty}\sigma(P_n,Z_n)=\frac{8-8\sqrt{2}+2\pi}{\pi}\approx 0.945.$$

If  $n \equiv 3 \pmod 4 $, then 
$$\lambda_k(Z_n)>\lambda_k(P_n),\ \ \  k=1,\ldots, \frac{n-3}{4},$$
 
 $$\lambda_k(Z_n)<\lambda_k(P_n),\ \ \ k=\frac{n+5}{4},\ldots, \frac{n-1}{2}, $$
 and $$\lambda_k(Z_n)=\lambda_k(P_n), \ \ \ k=\frac{n+1}{4},\frac{n+1}{2}.$$
 It follows from  Table 1 and the bipartiteness of $P_n$ and $Z_n$ that
  \begin{eqnarray*}
\sigma(P_n,Z_n)&=&2\big(\displaystyle\sum_{k=1}^{\frac{n-3}{4}}(\lambda_k(Z_n)-\lambda_{k}(P_n))+\displaystyle\sum_{k=\frac{n+5}{4}}^{\frac{n-1}{2}}(\lambda_k(P_n)-\lambda_{k}(Z_n))\big)\\
&=&4\big(\displaystyle\sum_{k=1}^{\frac{n-3}{4}}(\cos(\frac{(2k-1)\pi}{2n-2})-\cos(\frac{k\pi}{n+1}))
+\displaystyle\sum_{k=\frac{n+5}{4}}^{\frac{n-1}{2}}(\cos(\frac{k\pi}{n+1})-\cos(\frac{(2k-1)\pi}{2n-2}))\big).
\end{eqnarray*} 
By using Maple \cite{m},
  $$\displaystyle\lim_{n\longrightarrow \infty}\sigma(P_n,Z_n)=\frac{8-8\sqrt{2}+2\pi}{\pi}\approx 0.945.$$
If 
$$n^* = \left\{
 \begin{array}{ll}
 \frac{n}{4} &  n\equiv 0 \pmod 4,\\
 \\
 \frac{n-2}{4} &  n\equiv 2 \pmod 4,
 \end{array}
\right.$$
then
$$\lambda_k(Z_n)>\lambda_k(P_n),\ \ \  k=1,\ldots, n^*,$$
 and
 $$\lambda_k(Z_n)<\lambda_k(P_n),\ \ \ k=n^*+1,\ldots, \frac{n}{2}. $$

 Therefore
  \begin{eqnarray*}
\sigma(P_n,Z_n)&=&2\big(\displaystyle\sum_{k=1}^{n^*}(\lambda_k(Z_n)-\lambda_{k}(P_n))+\displaystyle\sum_{k=n^*+1}^{\frac{n}{2}}(\lambda_k(P_n)-\lambda_{k}(Z_n))\big )\\
&=&4\big(\displaystyle\sum_{k=1}^{n^*}(\cos(\frac{(2k-1)\pi}{2n-2})-\cos(\frac{k\pi}{n+1}))
+\displaystyle\sum_{k=n^*+1}^{\frac{n}{2}}(\cos(\frac{k\pi}{n+1})-\cos(\frac{(2k-1)\pi}{2n-2}))\big).
\end{eqnarray*}
 
In both cases, by using Maple \cite{m}, we have
 $$\displaystyle\lim_{n\longrightarrow \infty}\sigma(P_n,Z_n)=\frac{8-8\sqrt{2}+2\pi}{\pi}\approx 0.945.$$
  Therefore $a_n=\sigma(P_n,Z_n)$ is convergent to $\frac{8-8\sqrt{2}+2\pi}{\pi}$ .

\item Let $n\geq 6$ and $b_n:=\sigma(P_n,Z_n)$. We prove 
$$\displaystyle\lim_{n\longrightarrow \infty}b_{4n+1}=\displaystyle\lim_{n\longrightarrow \infty}b_{4n+2}=\displaystyle\lim_{n\longrightarrow \infty}b_{4n+3}=\displaystyle\lim_{n\longrightarrow \infty}b_{4n}=\frac{8-8\sqrt{2}+2\pi}{\pi},$$
which follows that $b_n$ is convergent to $\frac{8-8\sqrt{2}+2\pi}{\pi}$. \\

Suppose that $\lambda_1(W_n)\geq\cdots\geq\lambda_n(W_n)$ and $\lambda_1(Z_n)\geq\cdots\geq\lambda_n(Z_n)$ are the eigenvalues of 
$W_n$ and $Z_n$, respectively. By similar arguments given in part (2), if $n\equiv 1 \pmod 4$, then 

\begin{eqnarray*}
\sigma(W_n,Z_n)&=&2\big(\displaystyle\sum_{k=1}^{\frac{n-1}{4}}(\lambda_k(W_n)-\lambda_{k}(Z_n))+\displaystyle\sum_{k=\frac{n+3}{4}}^{\frac{n-1}{2}}(\lambda_k(Z_n)-\lambda_{k}(W_n))\big )\\
&=&4\big(\displaystyle\sum_{k=1}^{\frac{n-1}{4}}(\cos(\frac{(k-1)\pi}{n-3})-\cos(\frac{(2k-1)\pi}{2n-2}))
+\displaystyle\sum_{k=\frac{n+3}{4}}^{\frac{n-1}{2}}(\cos(\frac{(2k-1)\pi}{2n-2})-\cos(\frac{(k-1)\pi}{n-3}))\big).
\end{eqnarray*}
We have 
 \begin{eqnarray*}
\sigma(W_n,Z_n)&=&2\big(\displaystyle\sum_{k=1}^{\frac{n-3}{4}}(\lambda_k(W_n)-\lambda_k(Z_n))+\displaystyle\sum_{k=\frac{n+5}{4}}^{\frac{n-1}{2}}(\lambda_k(Z_n)-\lambda_k(W_n))\big)\\
&=&4\big(\displaystyle\sum_{k=1}^{\frac{n-1}{4}}(\cos(\frac{(k-1)\pi}{n-3})-\cos(\frac{(2k-1)\pi}{2n-2}))
+\displaystyle\sum_{k=\frac{n+5}{4}}^{\frac{n-1}{2}}(\cos(\frac{(2k-1)\pi}{2n-2})-\cos(\frac{(k-1)\pi}{n-3}))\big).
\end{eqnarray*} 
where $n \equiv 3 \pmod 4 $. In both cases, by using Maple \cite{m}, 
 $$\displaystyle\lim_{n\longrightarrow \infty}\sigma(W_n,Z_n)=\frac{8-8\sqrt{2}+2\pi}{\pi}\approx 0.945.$$
If
$$n^* = \left\{
 \begin{array}{ll}
 \frac{n}{4} &  n\equiv 0 \pmod 4,\\
 \\
 \frac{n-2}{4} &  n\equiv 2 \pmod 4,
 \end{array}
\right.$$
then

$$\lambda_k(W_n)>\lambda_k(Z_n),\ \ \  k=1,\ldots, n^*,$$
 
 $$\lambda_k(W_n)<\lambda_k(Z_n),\ \ \ k=n^*+1,\ldots, \frac{n}{2}-1, $$
 and
 $$\lambda_{\frac{n}{2}}(W_n)=\lambda_{\frac{n}{2}}(Z_n).$$ 

Therefore 
\begin{eqnarray*}
\sigma(W_n,Z_n)&=&2\big(\displaystyle\sum_{k=1}^{n^*}(\lambda_k(W_n)-\lambda_{k}(Z_n))+\displaystyle\sum_{k=n^*+1}^{\frac{n}{2}-1}(\lambda_k(Z_n)-\lambda_{k}(W_n))\big )\\
&=&4\big(\displaystyle\sum_{k=1}^{n^*}(\cos(\frac{(k-1)\pi}{n-3})-\cos(\frac{(2k-1)\pi}{2n-2}))
+\displaystyle\sum_{k=n^*+1}^{\frac{n}{2}-1}(\cos(\frac{(2k-1)\pi}{2n-2})-\cos(\frac{(k-1)\pi}{n-3}))\big).
\end{eqnarray*}
 
In both cases, by using Maple \cite{m},
 $$\displaystyle\lim_{n\longrightarrow \infty}\sigma(W_n,Z_n)=\frac{8-8\sqrt{2}+2\pi}{\pi}\approx 0.945.$$
Thus $b_n=\sigma(W_n,Z_n)$ is convergent to $\frac{8-8\sqrt{2}+2\pi}{\pi}$. 

\item The eigenvalues of $P_n$ and $W_n$ satisfy the same (in)equalities as the eigenvalues of $P_n$ and $Z_n$. So, by the parts (2) and (3), we obtain 
 
$$ \sigma(P_n,W_n)=\sigma(P_n,Z_n)+\sigma(W_n,Z_n).$$
Since 
$$\lim_{n\longrightarrow\infty}\sigma(P_n,Z_n)=\lim_{n\longrightarrow\infty}\sigma(W_n,Z_n),$$
we have
$$ \displaystyle\lim_{n\longrightarrow \infty}\sigma(P_n,Z_n)=\displaystyle\lim_{n\longrightarrow \infty}\sigma(W_n,Z_n)=\frac{1}{2}\displaystyle\lim_{n\longrightarrow \infty}\sigma(P_n,W_n)=\frac{8-8\sqrt{2}+2\pi}{\pi}\approx 0.945.$$

This completes the proof. $\hfill\Box$\\

\end{enumerate}

\section*{{\bf Acknowledgments}}

The research of the first author was in part supported by a grant from School of Mathematics, Institute for Research in Fundamental Sciences (IPM).

% ----------------------------------------------------------------

\end{document}